\documentclass{elsart}

\usepackage{natbib}

\usepackage{amssymb,graphicx}

\begin{document}

\begin{frontmatter}



\title{Modification of Tukey's Additivity Test}

\thanks{This research was supported by the grants MZE 0002701404 and NAZV QH81312. We are indebted to Dieter Rasch and an anonymous reviewer for their comments to the text.}

\author{Petr \v{S}ime\v{c}ek, Marie \v{S}ime\v{c}kov\'a}
\ead{simecek@gmail.com, simeckova.marie@vuzv.cz}


\address{Institute of Animal Science, P\v{r}\'{a}telstv\'{\i} 815, 10400 Prague, Czech Republic}

\begin{abstract}
In this paper we discuss testing for an interaction in the two--way ANOVA with just one observation per cell. The known results are reviewed and a simulation study is performed to evaluate type I and type II risks of the tests.
It is shown that the Tukey and Mandel additivity tests have very low power in case of more general interaction scheme. A modification of Tukey's test is developed to resolve this issue.
All tests mentioned in the paper have been implemented in R package AdditivityTests.
\end{abstract}


\end{frontmatter}

\section{Introduction}

In many applications of statistical methods, it is assumed that the response variable is a sum of several factors and a random noise. In a real world this may not be an appropriate model. For example, some patients may react differently to the same drug treatment or the influence of fertilizer may be influenced by the type of a soil. There might exist an interaction between factors. A testing for such interaction will be referred here as \textbf{testing of additivity hypothesis}.

If there is more than one observation per cell then standard ANOVA techniques may be applied. Unfortunately, in many cases it is infeasible to get more than one observation taken under the same conditions. For instance, it is not logical to ask the same student the same question twice. 
\medskip

We restrict ourselves to a case of two factors, i.e. two--array model, when the response in $i^{th}$ row and $j^{th}$ column is modeled as
\begin{equation}\label{m1}
y_{ij}=\mu+\alpha_i+\beta_j+\gamma_{ij}+\epsilon_{ij},\qquad i=1,\dots,a,\ j=1,\dots,b,
\end{equation}
where 
$$\sum_i \alpha_i = \sum_j \beta_j = \sum_i \gamma_{ij} = \sum_j \gamma_{ij} = 0$$
and the $\epsilon_{ij}$ are  normally distributed independent random variables with zero mean and variance $\sigma^2$.
\medskip

To test the additivity hypothesis
\begin{equation}\label{hyp}
  H_0\textrm{: }\gamma_{ij}=0\ \ i=1,\dots,a,\ j=1,\dots,b,
\end{equation}
a number of tests have been developed. The Section 2 recollects the known additivity tests, see also \citet{AlinKurt2006} and \citet{boik1993}.

In Section 3 the power of the tests described in Section 2 is compared by means of simulation. While Tukey test has relatively good power when the interaction is a product of the main effects, i.e. when $\gamma_{ij}=k \alpha_i \beta_j$ ($k$ is a real constant), its power for more general interaction is very poor. 

It should be reminded that \citet{tukey1949}  did not originally propose his
test for any particular type of interaction. Actually after a small
modification derived in Section 4 the power of the test improves dramatically.
There exist some issues when a sample size is not large enough that may be
resolve by a~permutation test or bootstrap.

\section{Additivity Tests}

This section recalls the known additivity tests of hypothesis~(\ref{hyp}) in
model~(\ref{m1}). Let $\bar{y}_{\cdot\cdot}$ denotes  the overall mean,
$\bar{y}_{i\cdot}$ the $i^{th}$ row's mean and $\bar{y}_{\cdot j}$ the $j^{th}$
column's mean. The matrix $R=[r_{ij}]$ will stand for a~residual matrix with respect to the main effects
$$r_{ij} = y_{ij} - \bar{y}_{i\cdot} - \bar{y}_{\cdot j} + \bar{y}_{\cdot\cdot}$$
The decreasingly ordered list of eigenvalues of $R R^T$ will be denoted by
$\kappa_1>\kappa_2>\dots>\kappa_{\min(a,b)-1}$, and its scaled versions equal  
$$\omega_i = \frac{\kappa_i}{\sum_k{\kappa_k}},\qquad i=1,2,\dots,\min(a,b)-1.$$
If the interaction is present we may expect that some of $\omega_i$ coefficients will be substantially higher than others.
\medskip

\noindent\textbf{Tukey test:} Introduced in \citet{tukey1949}. Tukey test first estimates row and column effects and then tests for the interaction of a type $\gamma_{ij}=k\alpha_i\beta_j$ ($k=0$ implies no interaction). Tukey test statistic $S_T$ equals 
$$S_T = \textsl{MS}_{int} / \textsl{MS}_{error},$$
where
$$\textsl{MS}_{int} = \frac{\left(\sum_i\sum_j y_{ij}(\bar{y}_{i\cdot}-\bar{y}_{\cdot\cdot})(\bar{y}_{\cdot j}-\bar{y}_{\cdot\cdot})\right)^2}{\sum_i(\bar{y}_{i\cdot}-\bar{y}_{\cdot\cdot})^2\sum_j(\bar{y}_{\cdot j}-\bar{y}_{\cdot\cdot})^2}$$ 
and
$$\textsl{MS}_{error} = \frac{\sum_i\sum_j(y_{ij}-\bar{y}_{\cdot\cdot})^2-a\sum_j(\bar{y}_{\cdot j}-\bar{y}_{\cdot\cdot})^2-b\sum_i(\bar{y}_{i\cdot}-\bar{y}_{\cdot\cdot})^2-\textsl{MS}_{int}}{(a-1)(b-1)-1}.$$
Under the additivity hypothesis $S_T$ is $F$-distributed with $1$ and $(a-1)(b-1)-1$ degrees of freedom. 
\medskip 

\noindent\textbf{Mandel test:} Introduced in \citet{Mandel1961}. Mandel test statistic $S_M$ equals
$$S_M=\frac{\sum_i(z_i-1)^2\sum_j(\bar{y}_{\cdot j}-\bar{y}_{\cdot\cdot})^2}{a-1}\ /\ \frac{\sum_i\sum_j\left((y_{ij}-\bar{y}_{i\cdot})-z_i(\bar{y}_{\cdot j}-\bar{y}_{\cdot\cdot})\right)^2}{(a-1)(b-2)},$$
where
$$z_i:=\frac{\sum_j y_{ij}(\bar{y}_{\cdot j}-\bar{y}_{\cdot\cdot})}{\sum_j (\bar{y}_{\cdot j}-\bar{y}_{\cdot\cdot})^2}.$$
Under the additivity hypothesis $S_M$ is $F$-distributed with  $a-1$ and $(a-1)\cdot(b-1)$ degrees of freedom. 
\medskip

Definitions of the three later tests slightly differ from their original versions. For $a$, $b$ fixed, a simulation may be used to get the critical values. 

\noindent\textbf{Johnson -- Graybill test:} Introduced in \citet{JohnsonGraybill1972}. Johnson -- Graybill test statistic is just
$S_{J}=\omega_1$. The additivity hypothesis is rejected if $S_{J}$ is high. 
\medskip

\noindent\textbf{Locally best invariant (LBI) test:} See \citet{Boik1993b}. LBI
test statistic equals (up to a monotonic transformation)
$$S_{L}=\sum_{i=1}^{\min(a,b)-1}{\omega}_i^2.$$
The additivity hypothesis is rejected if $S_{L}$ is high. 
\medskip

\noindent\textbf{Tusell test:} See \citet{Tusell1990}. Tusell test statistic
equals (up to a constant) 
$$S_{U}=\prod_{i=1}^{\min(a,b)-1}\omega_i.$$
The additivity hypothesis is rejected if $S_{U}$ is low. 
\medskip

As will be verified in the next section, Tukey and Mandel tests are appropriate if $\gamma_{ij}=k\alpha_i\beta_j$ while  Johnson -- Graybill, LBI and Tusell omnibus tests are suitable in cases of more complexed interactions.

\section{Simulation Study}

In this section simulation results about power of the additivity tests are presented. 
According to \citet{marie} the type-I-risk of the tests mentioned in Section~2 is not touched even when one of the effects in (\ref{m1}) is considered as random. 
The mixed effects model used for the simulation study is as (\ref{m1}) 
where $\mu$, $\alpha_i$ are constants, and $\beta_j$ are independent normally distributed random variables with zero mean and variance $\sigma_{\beta}^2$.

Two possible interaction schemes were under inspection:
\begin{itemize}
\item[A)] $\gamma_{ij} = k \alpha_i \beta_j$ where $k$ is a real constant.
\item[B)] $\gamma_{ij} = k \alpha_i \delta_j$ where $\delta_j$ are independent normally distributed random variables with zero mean and  variance $\sigma_{\beta}^2$, independent of $\beta_j$ and $\epsilon_{ij}$, and $k$ a real constant.
\end{itemize}
The $\epsilon_{ij}$ are independent normally distributed random variables with zero mean, $\mu=0$, and unit variance, $\sigma^2=1$.

The other parameters are equal to $\mu=0$, $\sigma_{\beta}^2=2$,  $\sigma^2=1$, $a=10$, 
$$(\alpha_1,\dots,\alpha_{10})=(-2.03,-1.92,-1.27,-0.70,0.46,0.61,0.84,0.94,1.07,2.00).$$
Two possibilities are considered for the $b$, either $b=10$ or $b=50$, and $10$ different values between $0$ and $12$ are considered for the interaction parameter~$k$.

For each combination of parameters' values a dataset was generated based on the model~(\ref{m1}), the tests of additivity were done and their results were noted down. The step was repeated $10\ 000$ times. The estimated power of the test is the percentage of the positive results. All tests were done on $\tilde\alpha=5\%$ level. 

The dependence of the power on~$k$ is visualized in Figure 1. As we can see,
while Tukey and Mandel tests outperformed omnibus tests for interaction A and
low $k$ and $b$, they completely fail to detect the interaction B even for a large value of $k$ and $b=50$. Therefore, it is desirable to develop a test which is able to detect a spectrum of practically relevant alternatives while still has the power comparable to the Tukey and Mandel tests for the most common interaction scheme A.

\begin{figure}
\begin{center}
\includegraphics[width=14cm]{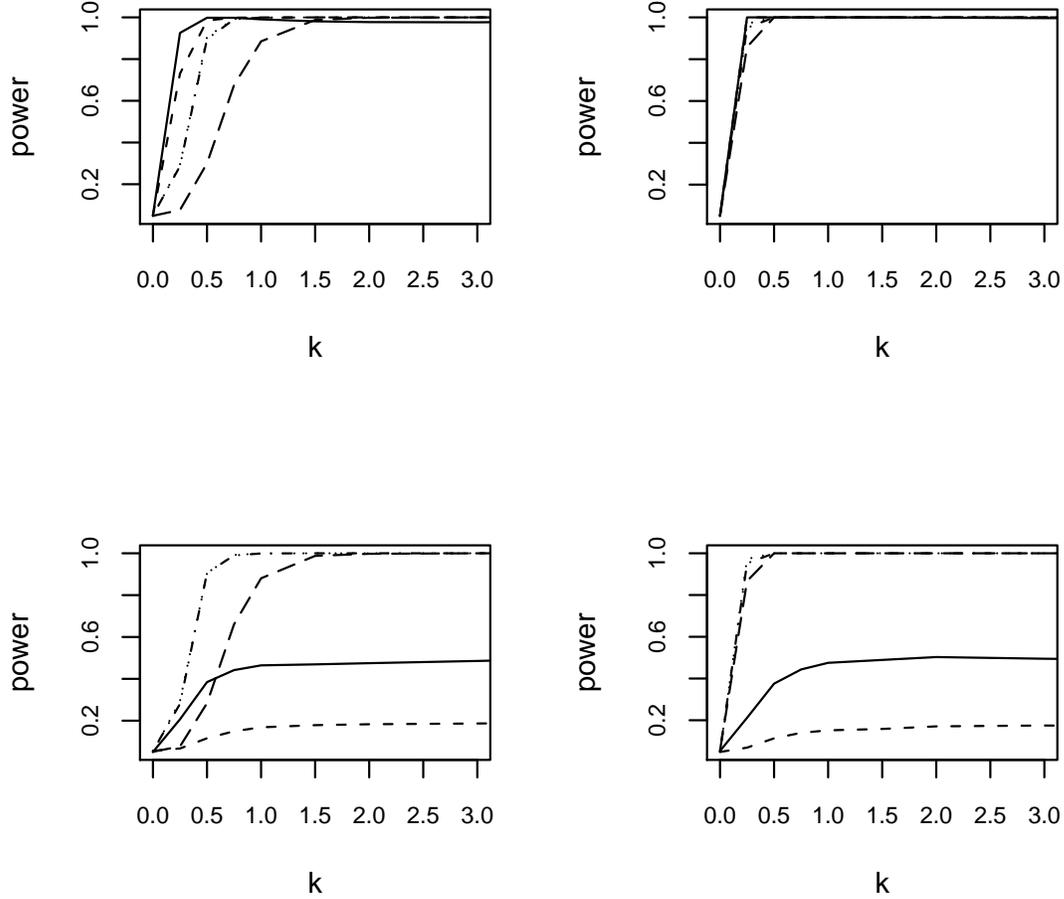}
\caption{Power dependence on $k$, $b$ ($b=10$ left, $b=50$ right) and 
interaction type ($A$ up, $B$ down). Tukey test solid line, Mandel test 
dashed line, Johnson -- Graybill test dotted line, LBI test dot-dash line, 
Tusell test long dash line.}
\end{center}
\end{figure}

\section{Modification of Tukey Test}

In Tukey test a model (\ref{m1})
\begin{equation}\label{modd}
y_{ij}=\mu+\alpha_i+\beta_j+\gamma_{ij}+\epsilon_{ij}=\mu+\alpha_i+\beta_j+k\alpha_i\beta_j+\epsilon_{ij}
\end{equation}
is tested against a submodel (\ref{hyp}) $y_{ij}=\mu+\alpha_i+\beta_j+\epsilon_{ij}$. The estimators of row effects $\hat{\alpha}_i=\bar{y}_{i\cdot}-\bar{y}_{\cdot\cdot}$ and column effects $\hat{\beta}_j=\bar{y}_{\cdot j}-\bar{y}_{\cdot\cdot}$ are calculated in the same way in both models although the dependency of $y_{ij}$ on these parameters is not linear for the full model.

The main idea behind a presented modification is that the full
model~(\ref{modd}) is fitted by a nonlinear regression and tested against a
submodel $y_{ij}=\mu+\alpha_i+\beta_j+\epsilon_{ij}$ by a likelihood ratio test. The estimates of row and column effects therefore differ for each model.

\subsection{Non-adjusted test}

Under additivity hypothesis the maximum likelihood estimators of parameters can be calculated simply as $\hat{\mu}=\bar{y}_{\cdot\cdot}$, 
$\hat{\alpha}_i=\bar{y}_{i\cdot}-\bar{y}_{\cdot\cdot}$ and $\hat{\beta}_j=\bar{y}_{\cdot j}-\bar{y}_{\cdot\cdot}$. Residual sum of squares equals
$$\textsl{RSS}_{0} = \sum_i\sum_j\left(y_{ij}-\hat{\mu}-\hat{\alpha}_i-\hat{\beta}_j\right)^2 = \sum_i\sum_j \left(y_{ij} - \bar{y}_{i\cdot} - \bar{y}_{\cdot j} + \bar{y}_{\cdot\cdot}\right)^2.$$

In the full model (\ref{modd}) the parameters' estimates are computed iteratively. Let us first take $\hat{\alpha}_i^{(0)}=\hat{\alpha}_i=\bar{y}_{i\cdot}-\bar{y}_{\cdot\cdot}$, $\hat{\beta}_j^{(0)}=\hat{\beta}_j=\bar{y}_{\cdot j}-\bar{y}_{\cdot\cdot}$ and $$\hat{k}^{(0)}=\frac{\sum_i\sum_j \left(y_{ij}-\hat{\alpha}^{(0)}_i-\hat{\beta}^{(0)}_j-\hat{\mu}\right)\cdot \hat{\alpha}^{(0)}_i\cdot \hat{\beta}^{(0)}_j}{\sum_i\sum_j \left(\hat{\alpha}^{(0)}_i\right)^2\cdot \left(\hat{\beta}^{(0)}_j\right)^2}.$$ The $\hat{k}^{(0)}$ is equal to the estimator of $k$ in the classical Tukey test.

The iteration procedure continues by updating estimates one by one (while the
rest of parameters are fixed):
\begin{itemize}
       \item $\hat{\alpha}_{i}^{(n)}=\frac{\sum_j \left(y_{ij}-\hat{\mu}-\hat{\beta}_j^{(n-1)}\right)\cdot \left(1+\hat{k}^{(n-1)}\cdot \hat{\beta}^{(n-1)}_j\right)}{\sum_j \left(1+\hat{k}^{(n-1)}\cdot \hat{\beta}^{(n-1)}_j\right)^2}$
       \item $\hat{\beta}_{j}^{(n)}=\frac{\sum_i \left(y_{ij}-\hat{\mu}-\hat{\alpha}^{(n-1)}_{i}\right)\cdot \left(1+\hat{k}^{(n-1)}\cdot \hat{\alpha}^{(n-1)}_{i}\right)}{\sum_i \left(1+\hat{k}^{(n-1)}\cdot \hat{\alpha}^{(n-1)}_{i}\right)^2}$
       \item $\hat{k}^{(n)}=\frac{\sum_i \sum_j \left(y_{ij}-\hat{\alpha}^{(n-1)}_{i}-\hat{\beta}^{(n-1)}_{j}-\hat{\mu}\right)\cdot \hat{\alpha}^{(n-1)}_{i}\cdot \hat{\beta}^{(n-1)}_{j}}
                                 {\sum_i \sum_j \left(\hat{\alpha}^{(n-1)}_{i}\right)^2\cdot \left(\hat{\beta}^{(n-1)}_{j}\right)^2}$
\end{itemize}
     
Surprisingly, it seems that one iteration is just enough to converge in a vast majority of cases. Therefore, for a simplicity reason let us define
$$\textsl{RSS} = \sum_i\sum_j\left(y_{ij}-\hat{\mu}-\hat{\alpha}_i^{(1)}-\hat{\beta}_j^{(1)}-k^{(1)}\hat{\alpha}_i^{(1)}\hat{\beta}_j^{(1)}\right)^2.$$
The likelihood ratio statistic of the modified Tukey test, i.e. a difference of twice log-likelihoods, equals
$$\frac{\textsl{RSS}_{0}-\textsl{RSS}}{\sigma^2}$$
and is asymptotically $\chi^2$-distributed with 1 degree of freedom. 

The consistent estimate of a residual variance $\sigma^2$ equals $s^2=\frac{\textsl{RSS}}{ab-a-b}$ and $\frac{\textsl{RSS}}{\sigma^2}$ is approximately $\chi^2$-distributed with $ab-a-b$ degrees of freedom. Thus, using a linear approximation of the nonlinear model (\ref{modd}) 
\begin{equation}\label{zmatek}
  \frac{\textsl{RSS}_{0}-\textsl{RSS}}{\frac{\textsl{RSS}}{ab-a-b}}
\end{equation}
is $F$-distributed with 1 and $ab-a-b$ degrees of freedom. 
Easy manipulation of (\ref{zmatek}) gives the modified Tukey test which rejects the additivity hypothesis if and only if 
$$\textsl{RSS}_{0}>\textsl{RSS} \left(1+\frac1{ab-a-b}F_{1,ab-a-b}(1-\tilde\alpha)\right),$$
where $F_{1,ab-a-b}(1-\tilde\alpha)$ stands for $1-\tilde\alpha$ quantile of $F$-distribution with $1$ and $ab-a-b$ degrees of freedom. 
\smallskip

\begin{figure}
\centering
\includegraphics[width=14cm]{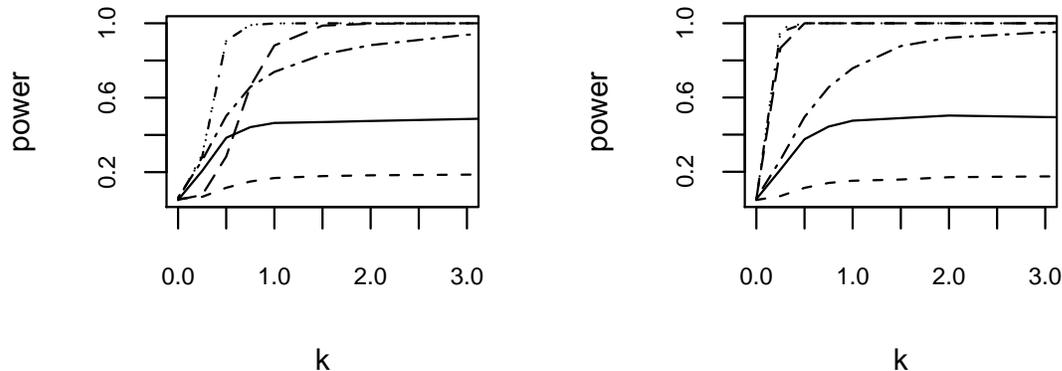}
\caption{Power dependence on $k$, $b$ ($b=10$ left, $b=50$ right) for 
interaction type $B$. Tukey test solid line, Mandel test dashed line, 
Johnson -- Graybill test dotted line, LBI test dot-dash line, Tusell test 
long dash line, modified Tukey test two dash line. The proposed modification
improved Tukey test and for large $k$ almost reach power of omnibus tests.}\label{obr2}
\end{figure}

Now we will return to the simulation study from Section~3. For interaction A the power of the modified test is almost equal to the power of Tukey test. For interaction B the power of the tests is compared on Figure~\ref{obr2}, the power of modified test is much higher than the power of Tukey test.

Theoretically, we may expect the modified test to be conservative because just one iteration does not find precisely the maximum of model~(\ref{modd}) likelihood. However, as we will see in the following part a situation for a small number of rows or columns is quite opposite.

\subsection{Small sample adjustment}

If the left part of Figure \ref{obr2} would be magnified enough it will show that the modified test does not work properly (type-I-risk $\doteq 6\%$). The reason is that the likelihood ratio test statistic converges to $\chi^2$-distribution rather slowly (see \citet{bartlet}) and a correction for small sample size is needed. We present two possibilities that are recommended if a number of rows or columns is below $20$ (empirical threshold based on simulations).
\medskip 

One possibility to overcome this obstacle is a permutation test, i.e. generate data as follows
$$y^{(perm)}_{ij}(t)=\hat{\mu}+\hat{\alpha}_i^{(0)}+\hat{\beta}_j^{(0)}+r_{\pi_{ij}(t)},\quad t=1,\dots,N^{(perm)}$$
where $\pi(t)$ is a random permutation of indexes of $R$ matrix. For each $t$ the statistic of interest $S^{(perm)}(t)=\textsl{RSS}_{0}(t)-\textsl{RSS}(t)$ is computed. The critical value equals 
$(1-\tilde{\alpha})\cdot 100\%$ quantile of $S^{(perm)}(t),\ t=1,\dots,N^{(perm)}$.
\medskip

The second possibility is to estimate the residual variance $s^2=\frac{\textsl{RSS}}{ab-a-b}$ and then generate samples of a distribution
$$y^{(sample)}_{ij}(t)=\hat{\mu}+\hat{\alpha}_i^{(0)}+\hat{\beta}_j^{(0)}+\epsilon^{(NEW)}_{ij}(t), \quad t=1,\dots,N^{(sample)}$$
where $(\epsilon^{(NEW)}_{ij})(t)$ are i.i.d. generated from a normal distribution with zero mean and variance $s^2$. This is simply parametric bootstrap on residuals.

The proposed statistic of interest is $\mathop{\rm abs}(k^{(1)})$ mirroring
deviation from null hypothesis $k=0$.
As in the permutation test the additivity hypothesis is rejected if more than $(1-\tilde{\alpha})\cdot 100\%$ of sampled statistics lie below the statistic based on real data. 

\section{Conclusion}

We have proposed a modification of the Tukey additivity test. The modified test performs almost as good as Tukey test when the interaction is a product of main effects but should be recommended if we also request  reasonable power in case of more general interaction schemes.
Problems with small sample size may be overcome by permutation test or parametric bootstrap on residuals.

All mentioned tests are implemented in R package \texttt{AdditivityTests} that may be downloaded on
$\texttt{http://github.com/rakosnicek/additivityTests}$.
As far as we are informed,  this is the first R implementation of additivity tests with the exception of the Tukey test.




\end{document}